\documentclass[aos,MSNbibl,seceqn,dvips]{arximspdf}

\doi{10.1214/13-AOS1136} 
\volume{41}
\issue{4}
\pubyear{2013}
\firstpage{2075}
\lastpage{2096}

\makeatletter

\newtheorem{theorem}{Theorem}[section]
\newtheorem{lemma}{Lemma}[section]
\newtheorem{corollary}{Corollary}[section]

\newproclaim{remark}{Remark}[section]

\newcommand{\E}{\mathrm{E}}
\newcommand{\R}{\mathbb{R}}

\newcommand{\boldX}{{\mathbf{X}}}
\newcommand{\boldR}{\mathbf{R}}
\newcommand{\boldY}{\mathbf{Y}}
\newcommand{\boldu}{\mu}
\newcommand{\boldv}{v}
\newcommand{\boldZ}{\mathbf{Z}}
\newcommand{\bolde}{e}
\newcommand{\tr}{\operatorname{tr}}
\newcommand{\idv}{\mathbf{1}}

\makeatother

\begin{document}
\begin{frontmatter}

\title{Tests for covariance matrix with fixed or divergent dimension}
\runtitle{Tests for HD covariance matrix}

\begin{aug}
\author[A]{\fnms{Rongmao} \snm{Zhang}\thanksref{T1}\ead[label=e3]{rmzhang@zju.edu.cn}},
\author[B]{\fnms{Liang} \snm{Peng}\thanksref{T2}\ead[label=e2]{peng@math.gatech.edu}}
\and
\author[C]{\fnms{Ruodu} \snm{Wang}\corref{}\thanksref{T3}\ead[label=e1]{wang@uwaterloo.ca}}
\runauthor{R. Zhang, L. Peng and R. Wang}
\affiliation{Zhejiang University, Georgia
Institute of Technology and University~of~Waterloo}
\address[A]{R. Zhang\\
Department of Mathematics\\
Zhejiang University\\
Hangzhou, Zhejiang 310027\\
China\\
\printead{e3}} 
\address[B]{L. Peng\\
School of Mathematics\\
Georgia Institute of Technology\\
Atlanta, Georgia 30332-0160\\
USA\\
\printead{e2}}
\address[C]{R. Wang\\
Department of Statistics and Actuarial Science\\
University of Waterloo\\
Waterloo, Ontario N2L 3G1\\
Canada\\
\printead{e1}}
\end{aug}

\thankstext{T1}{Supported by NSFC Grants 11171074 and
10801118, and the Fundamental Research Funds for the Central
Universities.}

\thankstext{T2}{Supported by NSF Grant DMS-10-05336.}

\thankstext{T3}{Supported in part by the Bob Price
Fellowship at the Georgia Institute of Technology.}

\received{\smonth{6} \syear{2012}}
\revised{\smonth{5} \syear{2013}}

%
\begin{abstract}
Testing covariance structure is of importance in many areas of
statistical analysis, such as microarray analysis and signal
processing. Conventional tests for finite-dimensional covariance
cannot be applied to high-dimensional data in general, and tests for
high-dimensional covariance in the literature usually depend on some
special structure of the matrix. In this paper, we propose some
empirical likelihood ratio tests for testing whether a covariance
matrix equals a given one or has a banded structure. The asymptotic
distributions of the new tests are independent of the dimension.
\end{abstract}

%

\end{frontmatter}

\section{Introduction}\label{sec1}

Let $\boldX_i=(X_{i1},\ldots, X_{ip})^T, i=1, 2,\ldots, n$, be
independent and identically distributed (i.i.d.) random vectors with
mean $\boldu=(\mu_1,\ldots, \mu_p)^T$ and covariance matrix
$\Sigma=(\sigma_{ij})_{1\le i,j\le p}$. For a given covariance matrix
$\Sigma_0$, it has been a long history for the study of testing
%
\begin{equation}
\label{11}H_0\dvtx  \Sigma=\Sigma_0 \quad\mbox{against}\quad
H_a\dvtx  \Sigma\neq\Sigma_0.
\end{equation}
Traditional methods for testing (\ref{11}) with finite $p$ include
the likelihood ratio test (see \cite{A2003}) and the scaled distance
measure for positive definite $\Sigma_0$ defined as
%
\begin{equation}
\label{12}V={1\over p}\operatorname{tr}(S_n-I_p
)^2,
\end{equation}
where $\operatorname{tr}(\cdot)$ denotes the trace of a matrix, $I_p$
denotes the $p\times p$ identity matrix and $S_n$ is the sample
covariance matrix of $\Sigma_0^{-1/2}\boldX_i$ (see \cite{J1971,J1972}
and~\cite{N1973}). When dealing with high-dimensional data,
the sample covariance in the likelihood ratio test is no longer
invertible with probability one, and tests based on a scaled
distance may also fail as demonstrated in \cite{LW2002}.

Since the above conventional tests cannot be employed for testing
high-dimensional covariance matrices, new methods are needed. When the
high-dimensional covariance matrix has a modest dimension $p$ compared
to the sample size $n$, that is, $p/n\rightarrow c$ for some
$c\in(0,\infty)$, Ledoit and Wolf \cite{LW2002} proposed a test by
modifying the scaled distance measure $V$ defined in (\ref{12}) under
the assumption that $\boldX_1$ has a normal distribution. When the
dimension $p$ is much larger than the sample size $n$, some special
structure has to be imposed. Chen, Zhang and Zhong \cite{CZZ2010a}
proposed a test which generalizes the result of \cite{LW2002} to the
case of nonnormal distribution and large dimension by assuming that
$\boldX_i=\Gamma\boldZ_i+\boldu$ for some i.i.d. $m$-dimensional random
vectors $\{\boldZ_i\}$ with $\mathrm{E}\boldZ_1=0$,
$\operatorname{var}(\boldZ_1)=I_m$, and $\Gamma$ is a $p\times m$
constant matrix with $\Gamma\Gamma^T=\Sigma$.

Sparsity is another commonly employed special structure in analyzing
high-dimensional data such as variable selection and covariance matrix
estimation.
Estimating sparse covariance
matrices has been actively studied in recent years. Some
recent references include \cite{BL2008,RLZ2009,CZZ2010} and \cite{CL2011}. When
the covariance matrix is assumed to be sparse and has a
banded structure, it becomes important to test whether the
covariance matrix possesses such a desired structure, that is, to test
%
\begin{equation}
\label{band} H_0\dvtx  \sigma_{ij}=0 \qquad\mbox{for all } |i-j|
\geq\tau,
\end{equation}
where $\tau<p$ is given and may depend on $n$.
Recently, Cai and Jiang \cite{CJ2011} proposed to use the maximum of the
absolute values of sample covariances to test (\ref{band}) when
$\boldX_1$ has a multivariate normal distribution. However, it is
known that the convergence rate of the normalized maximum to a
Gumbel limit is very slow, which means such a test has a poor
size in general. Although using maximum is very powerful in detecting
the departure from the null hypothesis when at least one large
departure exists, it is much less powerful than a test based on a
Euclidean distance when many small departures from the null hypothesis happen.

To avoid assuming the sparse structure and normality condition in the
testing problems (\ref{11}) and (\ref{band}), we propose to
construct tests based on the equivalent testing problem $H_0\dvtx
\|\Sigma-\Sigma_0\|^2_F =0$ against $H_a\dvtx
\|\Sigma-\Sigma_0\|^2_F\neq0$, where $\|A\|_F=\sqrt{\tr(A^TA)}$ is
the Frobenius norm of the matrix $A$.

Put $\boldY_i=(\boldX_i-\boldu)(\boldX_i-\boldu)^T$ for $i=1,\ldots, n$. Based on the fact
that $\E[\boldY_i]=\Sigma$, one can test (\ref{11}) by employing
the well-known Hotelling one-sample
$T^2$ statistic for a mean vector when $p$ is finite, and its modified
versions when
$p$ is divergent and some specific models are assumed; see,
for example, \cite{BS1996}~and~\cite{CQ2010}.

Another popular test for a finite-dimensional mean vector is the
empirical likelihood ratio test
proposed in \cite{O1988,O1990}. Recently, Hjort, McKeague and Van
Keilegom \cite{HMV2009} and Chen, Peng and Qin \cite{CPQ2009}
extended it to the high-dimensional case. It turns out that the asymptotic
distribution of the empirical likelihood ratio test is a chi-square
distribution for a fixed dimension and a normal distribution for a
divergent dimension. That is, the limit depends on
whether the dimension is fixed or divergent. Note that the methods in
the above papers can also be used to construct an estimator for unknown
parameters, which is called maximum empirical likelihood estimator.

Motivated by the empirical likelihood ratio test in \cite{PQW2013} for
testing a high-dimensional mean vector,
we propose to apply the empirical likelihood ratio test to
two estimating equations, where one equation ensures the consistency of
the proposed test and another one is used to improve the test power.
It turns out that the proposed test
puts no restriction on the sparse structure of the covariance
matrix and normality of $\boldX_1$. When testing (\ref{band}), a
similar procedure can be employed; see Section \ref{sec2} for more details.

The paper is organized as follows. In Section \ref{sec2}, we introduce the new
methodologies and present the main results. Simulation results are
given in Section \ref{sec3}. Section \ref{sec4} proves the main results. Detailed proofs
for lemmas used in Section \ref{sec4} are put in the supplementary material
\cite{ZPW2013}.

\section{Methodologies and main results}\label{sec2}

\subsection{Testing a covariance matrix}\label{sec2.1}\label{main}
Let $\boldX_i=(X_{i1},\ldots, X_{ip})^T, i=1,\ldots, n$, be
independent\vspace*{1pt} and identically distributed observations with mean $\boldu
=(\mu_1,\ldots, \mu_p)^T$ and covariance matrix
$\Sigma=(\sigma_{ij})_{p\times p}$.

When $\boldu$ is known, for $i=1,\ldots, n$, we define
$\boldY_i=(\boldX_i-\boldu)(\boldX_i-\boldu)^T$.
Then $\mathrm{E}[\tr((\boldY_1-\Sigma_0)(\boldY_2-\Sigma_0))]=0$
is equivalent to $\|\Sigma-\Sigma_0\|^2_F=0$, which is equivalent to
$H_0\dvtx  \Sigma=\Sigma_0$.
A direct application of the empirical likelihood ratio test to the above
estimating equation may endure low power by noting that
$\mathrm{E}[\tr((\boldY_1-\Sigma_0)(\boldY_2-\Sigma_0))]=\|\Sigma
-\Sigma_0\|^2_F=O(\delta^2)$ rather than $O(\delta)$ if
$\|\Sigma-\Sigma_0\|_F=O(\delta)$ and $p$ is fixed.
A brief simulation study and the power analysis in Section \ref{power}
confirm this fact. In order
to improve the test power, we propose to add one more linear
equation. Note that with prior information on the model or more
specific alternative hypothesis, a more proper linear equation may
be available. Without additional information, any linear equation
that detects the change of order $\|\Sigma-\Sigma_0\|_F$ is a
possible choice theoretically. Here we simply choose the following
functional ${\mathbf1}_{p}^T(\boldY_1+\boldY_2-2\Sigma_0)\idv_p$,
where ${\mathbf1}_{p}=(1,\ldots, 1)^T\in\R^{p}$. More
specifically, we propose to apply the empirical likelihood ratio test to
the following two equations:
%
\begin{equation}
\label{16}\qquad\mathrm{E}\bigl[\tr\bigl((\boldY_1-\Sigma_0)
(\boldY_2-\Sigma_0)\bigr)\bigr]=0\quad\mbox{and}\quad
\mathrm{E} \bigl[\idv_p^T(\boldY_1+
\boldY_2-2\Sigma_0)\idv_p\bigr]=0.
\end{equation}
Of course one can
try other linear equations or add more equations to further improve
the power. Theorems derived below can easily be extended to the case
when $\idv_p$ is replaced by any constant vector.

In order to obtain an independent paired data $(\boldY_1,
\boldY_2)$, we split the sample into two subsamples with size
$N=[n/2]$. That is, for $i=1, 2,\ldots, N$, we define
$\boldR_i(\Sigma)=(\bolde_i(\Sigma), \boldv_i(\Sigma))^T$, where
\[
\bolde_i(\Sigma)=\tr\bigl((\boldY_i-\Sigma) (
\boldY_{i+N}-\Sigma)\bigr)\quad\mbox{and}\quad \boldv_i(\Sigma)={
\mathbf1}_{p}^T(\boldY_i+\boldY_{i+N}-2
\Sigma)\idv_p.
\]
Based on $\{\boldR_i(\Sigma)\}_{i=1}^{N}$, we define the empirical
likelihood ratio function for $\Sigma$ as
\[
L_1(\Sigma)=\sup\Biggl\{\prod
_{i=1}^{N}(Np_i)\dvtx  \sum
_{i=1}^{N}p_i=1, \sum
_{i=1}^{N} p_i\boldR_{i}(
\Sigma)=0, p_1\geq0,\ldots,p_N\ge0\Biggr\}.
\]

When $\boldu$ is unknown, instead of using $\{\boldR_i(\Sigma)\}
_{i=1}^N$, we use $\{\boldR^*_i(\Sigma)\}_{i=1}^N$ where $\boldu$ is
replaced by the sample means.
That is, put $\overline{\boldX}{}^1={1\over N}\sum_{i=1}^{N}
\boldX_{i}$, $\overline{\boldX}{}^2={1\over N}\sum_{i=N+1}^{2N}
\boldX_{i}$,
and define
\[
\boldY^*_{i} =\bigl(\boldX_i-\overline{
\boldX}{}^1\bigr) \bigl(\boldX_i-\overline{
\boldX}{}^1\bigr)^T \quad\mbox{and}\quad \boldY^*_{N+i}=
\bigl(\boldX_{N+i}-\overline{\boldX}{}^2\bigr) \bigl(
\boldX_{N+i}-\overline{\boldX}{}^2\bigr)^T
\]
for $i=1,\ldots,N$. Put $\boldR^*_i(\Sigma)=(\bolde^*_i(\Sigma),
\boldv^*_i(\Sigma))^T$, where
\[
\bolde^*_i(\Sigma)=\tr\biggl(\biggl(\boldY_i^*-
{(N-1)\Sigma\over N}\biggr) \biggl(\boldY^*_{i+N}-
{(N-1)\Sigma\over N}\biggr)\biggr)
\]
and
\[
\boldv^*_i(\Sigma)={\mathbf1}_{p}^T\biggl(
\boldY^*_i+\boldY^*_{i+N}-{2(N-1)\Sigma\over N}\biggr)
\idv_p.
\]
As before, we define the empirical likelihood ratio function for
$\Sigma$ as
\[
L_2(\Sigma)=\sup\Biggl\{\prod_{i=1}^{N}(Np_i)\dvtx
\sum_{i=1}^{N}p_i=1, \sum
_{i=1}^{N} p_i
\boldR^*_{i}(\Sigma)=0, p_1\geq0,\ldots,p_N\ge
0\Biggr\}.
\]

First we show that Wilks's theorem holds for the above empirical
likelihood ratio tests without imposing any special structure on
$\boldX_1$.

\begin{theorem} \label{main1} Suppose that $\mathrm{E}[\boldv
_1^2(\Sigma)]>0$ and for some $\delta>0$,
%
\begin{eqnarray}
\label{t21}\quad&&\max\bigl\{{\mathrm{E}\bigl|\bolde_1(\Sigma)\bigr|^{2+\delta}/
\bigl(\E\bigl[e_1^2(\Sigma)\bigr]\bigr)^{({2+\delta})/2}},
{\mathrm{E}|\boldv_1(\Sigma)|^{2+\delta}/ \bigl(\E
\bigl[v_1^2(\Sigma)\bigr]\bigr)^{({2+\delta})/
2}} \bigr\}
\nonumber\\[-8pt]\\[-8pt]
&&\qquad=o\bigl(N^{({\delta+\min\{2, \delta\}})/4}\bigr).
\nonumber
\end{eqnarray}
Then under $H_0\dvtx  \Sigma=\Sigma_0$, $-2\log L_1(\Sigma_0)$
converges in distribution to a chi-square distribution with two
degrees of freedom as $n\rightarrow\infty$.
In addition, if
%
\begin{equation}
\label{t21add}\qquad\bigl(\tr\bigl(\Sigma^2\bigr)\bigr)^2=o
\bigl(N^2\mathrm{E}\bigl[e_1^2(\Sigma)\bigr]
\bigr)  \quad\mbox{and}\quad \Biggl(\sum_{i=1}^p\sum
_{j=1}^{p}\sigma_{ij}
\Biggr)^2=o\bigl(N\mathrm{E}\bigl[v_1^2(
\Sigma)\bigr]\bigr),
\end{equation}
then under $H_0\dvtx  \Sigma=\Sigma_0$, $-2\log L_2(\Sigma_0)$
also converges in distribution to a chi-square distribution with two
degrees of freedom as $n\rightarrow\infty$.
\end{theorem}

Using Theorem \ref{main1}, one can test $H_0\dvtx  \Sigma=\Sigma_0$
against $H_a\dvtx  \Sigma\neq\Sigma_0$. That is, one rejects $H_0$ at
level $\alpha$ when $-2\log L_1(\Sigma_0)>\xi_{1-\alpha}$ if $\mu$
is known, or when $-2\log L_2(\Sigma_0)>\xi_{1-\alpha}$ if $\mu$ is
unknown, where
$\xi_{1-\alpha}$ denotes the $(1-\alpha)$th quantile of a chi-square
distribution with two
degrees of freedom.

Write the $p\times p$ matrix $\boldY_1$ as a $q=p^2$-dimensional
vector, and denote the covariance matrix of such a vector by
$\Theta=(\theta_{ij})_{q\times q}$. Conditions in Theorem~\ref
{main1} can
be guaranteed by imposing some conditions on the moments and
dimensionality of $\boldX_1$ such as the following assumptions:
\begin{itemize}
\item[A1:] $\liminf_{n\rightarrow\infty}\frac1q\operatorname{tr}(\Theta
^2)>0$ and $\liminf_{n\rightarrow\infty}\frac1q\idv_q^T\Theta\idv_q>0$;

\item[A2:] For some $\delta>0$, ${1\over p^2}\sum_{i=1}^{p}\sum
_{j=1}^{p}\mathrm{E}|(X_{1, i}-\mu_i)(X_{1, j}-\mu_j)-\sigma
_{ij}|^{2+\delta}=O(1)$;

\item[A3:] $p=o (n^{({\delta+\min(2, \delta)})/({4(2+\delta
)})} )$.
\end{itemize}

\begin{corollary} \label{main2} Under conditions \textup{A1--A3} and $H_0\dvtx
\Sigma=\Sigma_0$,\break $-2\log L_1(\Sigma_0)$ converges in distribution to
a chi-square distribution with two degrees of freedom as
$n\rightarrow\infty$.
Further, if
%
\begin{equation}
\label{th2add1} \max_{1\leq i\leq p}\sigma_{ii}<C_0
\qquad\mbox{for some constant } C_0>0,
\end{equation}
then $-2\log L_2(\Sigma_0)$ also converges in distribution to
a chi-square distribution with two degrees of freedom as
$n\rightarrow\infty$.
\end{corollary}

\begin{remark}
\label{rem1}
Condition (\ref{t21}) requires that the second moment of $(e_1,v_1)$
is not too small compared to a higher-order moment of $(e_1,v_1)$,
which ensures
that Lyapunov central limit theorem holds for ${1\over
\sqrt{N}}\sum_{i=1}^{N}\bolde_i(\Sigma_0)$ and ${1\over
\sqrt{N}}\sum_{i=1}^{N}\boldv_i(\Sigma_0)$. Condition (\ref
{t21add}) makes sure that the mean vector can be replaced by the
sample mean.
It is obvious that (\ref{t21add}) and (\ref{th2add1})
hold when $p$ is fixed.
\end{remark}

Note that condition A1 is only related to the covariance matrix and
condition~A2 holds obviously if
\[
{1\over p^2}\sum_{i=1}^{p}\sum
_{j=1}^{p}\mathrm{E}|X_{1,
i}X_{1,j}|^{2+\delta}<
\infty\quad\mbox{or}\quad {1\over p} \sum_{i=1}^{p}
\mathrm{E}|X_{1i}|^{4+2\delta}<\infty.
\]
Condition A3 imposes some restriction on $p$, but it can be removed if
$\boldX_i$ has some special dependence structure. For example, Theorem
\ref{main1} can be applied to the following setting studied in \cite{CQ2010,BS1996} and \cite{CZZ2010a}:

(B) (Multivariate model). Assume that the sample has
the following decomposition:
%
\begin{equation}
\label{29}\boldX_i=\Gamma\boldZ_i+\boldu,
\end{equation}
where
$\Gamma$ is a $p\times m$ constant matrix with $\Gamma\Gamma
^T=\Sigma$ and $\{\boldZ_i=(Z_{i1},\ldots,Z_{im})^T\}$ is a sequence
of $m$-dimensional i.i.d. random vectors with $\E\boldZ_i=0,
\operatorname{var}(\boldZ_i)=I_m, \mathrm{E}(Z_{11}^4 ) =\cdots
=\mathrm{E}(Z_{1m}^4)=3+\Delta>1$ and uniformly bounded 8th moment.
Further assume that for any integers $l_v\ge0$ and $h\ge1$ with
$\sum_{v=1}^{h}l_v=8$,
%
\begin{equation}\label{moments}
\mathrm{E}\bigl(Z_{1i_1}^{l_1}Z_{1i_2}^{l_2}
\cdots Z_{1i_{h}}^{l_h}\bigr)=\mathrm{E}\bigl(Z_{1i_1}^{l_1}
\bigr) \mathrm{E}\bigl(Z_{1i_2}^{l_2}\bigr)\cdots\mathrm{E}
\bigl(Z_{1i_h}^{l_h}\bigr),
\end{equation}
where $i_1,\ldots,i_h$ are distinct.

Note that if $\boldX_i$ has a multivariate normal distribution, then
(B) holds.

\begin{corollary} \label{main3} Suppose \textup{(B)} holds with $\sum
_{i=1}^{p}\sum_{j=1}^{p}\sigma_{ij}>0$. Then, under $H_0\dvtx  \Sigma
=\Sigma_0$, both $-2\log L_1(\Sigma_0)$ and $-2\log L_2(\Sigma_0)$
converge in distribution to
a chi-square distribution with two degrees of freedom as
$n\to\infty$.
\end{corollary}

\begin{remark}\label{rem2}
Note that condition $\sum_{i=1}^p\sum_{j=1}^p\sigma_{ij}>0$ for
model (B) implies that $\mathrm{E}[\boldv_1^2(\Sigma)]>0$; see
the proof of Lemma \ref{lem3}.
For testing $H_0\dvtx \Sigma=I_p$, \cite{CZZ2010} proposed a test based
on the above model and required $p\to\infty$ as $n\to\infty$. In
comparison, the proposed empirical likelihood ratio tests work for both
fixed and divergent~$p$.
\end{remark}
%
\begin{remark}
\label{rem3}
When one is interested in testing $H_0\dvtx  \mu=\mu_0$ and $\Sigma=\Sigma
_0$, it is straightforward to combine the proposed empirical likelihood
ratio test with that in \cite{PQW2013} for testing a high-dimensional mean.
\end{remark}

\subsection{Testing bandedness}\label{sec2.2} Suppose $\{\boldX_i\}$ is a sequence
of i.i.d. normal random vectors with covariance matrix $\Sigma=(\sigma
_{ij})_{1\le i,j\le p}$. Cai and Jiang \cite{CJ2011} proposed to use
the maximum of the absolute values of the
sample correlations to test a banded
structure
%
\begin{equation}
\label{210}H_0\dvtx  \sigma_{ij}=0 \qquad\mbox{for all } |i-j|\geq
\tau,
\end{equation}
where $\tau<p$. It is known that
the rate of convergence of the above maximum to a Gumbel distribution is
very slow in general, which results in a poor size; see also the
simulation results in Section \ref{sim}. Using the maximum as a test
statistic is powerful when at least a large deviation from the null
hypothesis exists. However, when many small deviations from the null
hypothesis exist, a test based on the maximum is much less efficient
than a test based on a Euclidean distance such as the test in \cite
{QC2012}. Here we modify
the empirical likelihood ratio tests in Section \ref{main} to test the
above banded
structure as follows.

For a matrix $M$, define the matrix $M^{(\tau)}$ as
$(M^{(\tau)})_{ij} = (M)_{ij} I(|i-j|{\geq}\tau)$, where $I(\cdot)$
denotes the indicator function.
Put
\begin{eqnarray*}
\bolde^{\prime}_i(\Sigma) &=& \tr\bigl(\bigl(
\boldY^{(\tau)}_i-\Sigma^{(\tau)}\bigr) \bigl(
\boldY^{(\tau
)}_{N+i}-\Sigma^{(\tau)}\bigr)\bigr), %
\\
{\boldv}^{\prime}_i(\Sigma)&=&{\mathbf1}_{p}^T
\bigl(\boldY^{(\tau)}_i+\boldY^{(\tau)}_{N+i}-2
\Sigma^{(\tau)}\bigr)\idv_p,
\\
\bolde^{*\prime}_i(\Sigma) &=& \tr\biggl(\biggl(
\boldY^{*(\tau)}_i-\frac{N-1}{N}\Sigma^{(\tau)}\biggr)
\biggl(\boldY^{*(\tau)}_{N+i} -\frac{N-1}{N}
\Sigma^{(\tau)}\biggr)\biggr)
\end{eqnarray*}
and
\[
{\boldv}^{*\prime}_i(\Sigma)={\mathbf1}_{p}^T
\biggl(\boldY^{*(\tau
)}_i+\boldY^{*(\tau)}_{N+i}-
\frac{2(N-1)}{N}\Sigma^{(\tau)}\biggr)\idv_p.
\]
Then $\Sigma^{(\tau)}$ is zero under $H_0$ in (\ref{210}).
Based on $\boldR^{\prime}_i(\Sigma)=(\bolde^{\prime}_i(\Sigma),
\boldv^{\prime}_i(\Sigma))^T$ and $\boldR^{*\prime}_i(\Sigma
)=(\bolde^{*\prime}_i(\Sigma), \boldv^{*\prime}_i(\Sigma))^T$,
we define the empirical likelihood ratio functions for $\Sigma$ as
\[
L_3(\Sigma)=\sup\Biggl\{\prod_{i=1}^{N}(Np_i)\dvtx
\sum_{i=1}^{N}p_i=1, \sum
_{i=1}^{N} p_i
\boldR^{\prime}_{i}(\Sigma)=0, p_i\geq0, i=1,\ldots, N\Biggr\}
\]
for the case of a known mean
and
\[
L_4(\Sigma)=\sup\Biggl\{\prod_{i=1}^{N}(Np_i)\dvtx
\sum_{i=1}^{N}p_i=1, \sum
_{i=1}^{N} p_i
\boldR^{*\prime}_{i}(\Sigma)=0, p_i\geq0, i=1,\ldots, N\Biggr\}
\]
for the case of an unknown mean.
Similar to the proof of Theorem \ref{main1}, we can show that
$-2\log L_3(\Sigma_0)$ and $-2\log L_4(\Sigma_0)$
converge in distribution to a chi-square distribution with two
degrees of freedom as $n\rightarrow\infty$ under some moment conditions.

\begin{theorem} \label{banded} Suppose that $\mathrm{E}[\boldv
'^{2}_1(\Sigma)]>0$ and for some $\delta>0$,
%
\begin{eqnarray}
\label{29c}\quad&&\max\bigl\{{\mathrm{E}\bigl|\bolde'_1(
\Sigma)\bigr|^{2+\delta}/ \bigl(\E\bigl[e_1'^2(
\Sigma)\bigr]\bigr)^{({2+\delta})/2}}, {\mathrm{E}\bigl|\boldv'_1(
\Sigma)\bigr|^{2+\delta}/ \bigl(\E\bigl[v_1'^2(
\Sigma)\bigr]\bigr)^{({2+\delta})/2}}\bigr\}
\nonumber\\[-8pt]\\[-8pt]
&&\qquad=o\bigl(N^{({\delta+\min\{2, \delta\}})/4}\bigr).
\nonumber
\end{eqnarray}
Then under $H_0$ in (\ref{210}), $-2\log L_3(\Sigma_0)$
converges in distribution to a chi-square distribution with two
degrees of freedom as $n\rightarrow\infty$, where $\Sigma_0$ is any
matrix such that $\Sigma_0^{(\tau)}=0$. In addition, if
\[
\E\Biggl\{ \sum_{i=1}^{N}
\bigl(e_i^{*\prime}(\Sigma)-e_i^{\prime}(
\Sigma)\bigr)^2+ \Biggl[\sum_{i=1}^{N}
\bigl(e_i^{*\prime}(\Sigma)-e_i^{\prime}(
\Sigma)\bigr)\Biggr]^2\Biggr\} =o\bigl(N\mathrm{E}
\bigl[e_1^{\prime2}(\Sigma)\bigr]\bigr)
\]
and
\[
\E\Biggl\{\sum_{i=1}^{N}
\bigl(v_i^{*\prime}(\Sigma)-v_i^{\prime}(
\Sigma)\bigr)^2+ \Biggl[\sum_{i=1}^{N}
\bigl(v_i^{*\prime}(\Sigma)-v_i^{\prime}(
\Sigma)\bigr)\Biggr]^2\Biggr\}= o\bigl(N\mathrm{E}
\bigl[v_1^{\prime2}(\Sigma)\bigr]\bigr),
\]
then under $H_0$ in (\ref{210}), $-2\log L_4(\Sigma_0)$ also
converges in
distribution to a chi-square distribution with two degrees of
freedom as $n\rightarrow\infty$.
\end{theorem}

In order to compare with \cite{CJ2011}, we use a different linear
functional so as to easily verify conditions when $\boldX_i$ has a
multivariate normal distribution. More specifically,
for a $p\times p$ matrix $M$, we define the matrix $M^{[\tau]}$ as
\begin{eqnarray*}
\bigl(M^{[\tau]}\bigr)_{ij}&=&(M)_{ij} \bigl\{I\bigl(i
\le(p-\tau)/2, j> (p+\tau)/2\bigr)\\
&&\hspace*{30pt}{}+I\bigl(j\le(p-\tau)/2, i> (p+\tau
)/2\bigr)
\bigr\}.
\end{eqnarray*}
Put
$
\tilde\boldv^{\prime}_i(\Sigma)={\mathbf1}_{p}^T(\boldY^{[\tau
]}_i+\boldY^{[\tau]}_{N+i}-2\Sigma^{[\tau]})\idv_p$
and
\[
\tilde\boldv^{*\prime}_i(\Sigma)={\mathbf1}_{p}^T
\biggl(\boldY^{*[\tau
]}_i+\boldY^{*[\tau]}_{N+i}-
{2(N-1)\over N}\Sigma^{[\tau]}\biggr)\idv_p.
\]
Based on
$\tilde\boldR^{*\prime}_i(\Sigma)=(\bolde^{*\prime}_i(\Sigma),
\tilde\boldv^{*\prime}_i(\Sigma))^T$, we define the empirical
likelihood ratio
function for $\Sigma$ as
\[
L_5(\Sigma)=\sup\Biggl\{\prod_{i=1}^{N}(Np_i)\dvtx
\sum_{i=1}^{N}p_i=1, \sum
_{i=1}^{N} p_i\tilde
\boldR^{*\prime}_{i}(\Sigma)=0, p_i\geq0, i=1,\ldots, N\Biggr\}.
\]

\begin{theorem}\label{THCJb} Assume $\boldX_i\sim N(\mu,\Sigma)$,
\[
C_1\le\liminf_{n\to\infty}\min_{1\le i\le p}
\sigma_{ii} \leq\limsup_{n\to\infty}\max
_{1\le i\le p}\sigma_{ii}\le C_2
\]
for some constants $0<C_1\le C_2<\infty$ and
$\tau=o( (\sum_{1\leq i, j\leq p}\sigma_{ij}^2)^{1/2})$. Then under
$H_0$ in (\ref{210}),
$-2\log L_5(\Sigma_0)$ converges in distribution to a chi-square
distribution with two degrees of freedom as
$n\to\infty$, where $\Sigma_0$ is any matrix such that $\Sigma
_0^{(\tau)}=0$.
\end{theorem}

\begin{remark}\label{rem4} Condition (\ref{29c}) is similar to (\ref
{t21}) to ensure that central limit theorem can be employed. The other
two conditions in Theorem \ref{banded} are similar to (\ref
{t21add}), and they make sure that the mean vector can be replaced by
the sample mean. The test in \cite{CJ2011}
requires that $\tau=o(p^{s})$ for all $s>0$ and $\log{p}=o(n^{1/3})$.
However, the new test in Theorem \ref{THCJb} only imposes conditions
between $\tau$ and
$p$. Also note that $\tau=o(p^{1/2})$ is sufficient for $\tau=o(
(\sum_{1\leq i, j\leq p}\sigma_{ij}^2)^{1/2})$.
\end{remark}

\subsection{Power analysis}\label{sec2.3}\label{power}

In this subsection we analyze the powers of our new tests. Denote
$\pi_{11}=\E(e_1^2(\Sigma))$, $\pi_{22}=\E(v_1^2(\Sigma))$,
$\zeta_{n1}=\operatorname{tr}((\Sigma-\Sigma_0)^2)/ \sqrt{\pi_{11}}$,
$\zeta_{n2}={2{\mathbf1}_{p}^T(\Sigma-\Sigma_0){\mathbf1}_p/
\sqrt{\pi_{22}}}$ and
$\nu=N(\zeta_{n1}^2+\zeta_{n2}^2)$.
Let $\xi_{\beta}$ denote the $\beta$-quantile of a chi-square
distribution with two degrees of freedom, and let $\chi_{2, \nu}^2$
denote a noncentral chi-square
distribution with two degrees of freedom and noncentrality parameter
$\nu$.

\begin{theorem}\label{power1} Under conditions of Corollary \ref
{main3} and $H_a\dvtx  \Sigma\neq\Sigma_0$, we have as $n\to\infty$,
%
\begin{equation}
\label{022} P\bigl\{-2\log L_j(\Sigma_0)>
\xi_{1-\alpha}\bigr\}=P\bigl\{\chi_{2,
\nu}^2>
\xi_{1-\alpha}\bigr\}+o(1)
\end{equation}
for $j=1, 2$.
\end{theorem}

\begin{remark}\label{rem5}
Note that under model (B), ${\pi_{11}}=O(\operatorname{tr}(\Sigma
^2)^2)$ and ${\pi_{22}}=O(\idv_p^T\Sigma\idv_p)^2$; see the proof
of Lemma \ref{lem3}. Therefore, $\zeta_{n1}=O(\operatorname{tr}((\Sigma-\Sigma
_0)^2)/\break\operatorname{tr}(\Sigma^2))$
and $\zeta_{n2}=O({\mathbf1}_{p}^T(\Sigma-\Sigma_0){\mathbf1}_p/
(\idv_p^T\Sigma\idv_p))$ are both natural measures of distance
between the null hypothesis and the real model.
\end{remark}

\begin{remark}
For a test only using the first estimating equation in (\ref{16}),
one needs $\sqrt{n}\zeta_{n1}\rightarrow\infty$ to ensure the
probability of rejecting $H_0$ goes to one. Thus it is less powerful
than the test using both estimating equations in (\ref{16}) when
$\sqrt{n}\zeta_{n2}\rightarrow\infty$ and $\sqrt{n}\zeta_{n1}$ is
bounded from infinity.
\end{remark}

From the above theorem, we conclude that the new test rejects $H_0$
with probability tending to one when either
$\sqrt{n}\zeta_{n1}$ or $\sqrt n|\zeta_{n2}|$ goes to infinity. To
compare with the power of the test given in \cite{CZZ2010a}, we
consider the testing problem $H_0\dvtx \Sigma=I_p$ against
$H_a\dvtx \Sigma\neq I_p$, where $\Sigma=I_p+(d I(|i-j|\leq\tau))_{1\le
i,j\le p}$ for some positive $d=d(n)\to0$ as $n\to\infty$. Note that
the term
$\sqrt{n^2\rho_{2,n}^2+n\rho_{2,n}}$ in (3.6) of \cite{CZZ2010a} is
a typo, and it should be $\sqrt{\rho_{2,n}^2+\rho_{2,n}}$. It is easy
to verify that the power of the test in \cite{CZZ2010a}
tends to one when $nd^2\tau\to\infty$ for the above example. On the
other hand, similar to Theorem 4 in \cite{CZZ2010a},
$\sqrt n|\zeta_{n2}|\rightarrow\infty$ is equivalent to $\sqrt
n|{2{\mathbf1}_{p}^T(\Sigma-\Sigma_0){\mathbf1}_p}|/p\rightarrow\infty$.
Thus the proposed empirical likelihood ratio test only needs
$nd^2\tau^2\to\infty$ to ensure that the power tends to one. Hence,
when $\Sigma=I_p+(d 1(|i-j|\le\tau))_{1\le i,j\le p}$ and
$\tau=\tau(n)\to\infty$, the proposed empirical likelihood ratio
test has
a larger local power than the test in \cite{CZZ2010a}.
For some other settings, the test in \cite{CZZ2010a} may be more powerful.

For testing\vspace*{1pt} the banded structure in Theorems \ref{banded} and
\ref{THCJb}, we have similar power results. Here we focus on Theorem
\ref{THCJb}. Let
$\kappa_{11}=\E(\bolde_1'^2(\Sigma))$ and
$\kappa_{22}=\E(\tilde\boldv_1'^2(\Sigma))$. Define
$\zeta'_{n1}=\operatorname{tr}((\Sigma^{(\tau)})^2)/
\sqrt{\kappa_{11}}$, $\zeta'_{n2}={2{\mathbf
1}_{{p}}^T\Sigma^{[\tau]}{\mathbf1}_p/ \sqrt{\kappa_{22}}}$ and $
\nu'=N(\zeta_{n1}'^2+\zeta_{n2}'^2)$.

\begin{theorem}\label{power2}
Under conditions of Theorem \ref{THCJb}, when $H_0$ in (\ref{210})
is false, we have as
$n\to\infty$
%
\begin{equation}
P\bigl\{-2\log L_5(\Sigma_0)>\xi_{1-\alpha}\bigr
\}=P\bigl\{\chi_{2,
\nu'}^2>\xi_{1-\alpha}\bigr\}+o(1),
\end{equation}
where $\Sigma_0$ is any matrix
such that
$\Sigma_0^{(\tau)}=0$.
\end{theorem}

\begin{remark}
As we argue in the \hyperref[sec1]{Introduction}, the size of the test in \cite{CJ2011}
is poor for testing a banded structure. Since the power
analysis for the test in \cite{CJ2011} is not available, theoretical
comparison is impossible. Instead, a simulation
comparison is given in the next section, which clearly shows that
the proposed test is much more powerful than the test in \cite{CJ2011}
when many small deviations from the null hypothesis exist. On the other
hand, the test in \cite{CJ2011} is more powerful when only a large
deviation exists. In that case, one can add more equations or replace
the second equation by a more relevant one in the proposed empirical
likelihood ratio test so as to catch this sparsity effectively.
\end{remark}

\section{Simulation}\label{sec3}\label{sim}
In this section we investigate the finite sample behavior of the
proposed empirical likelihood ratio tests in terms of both size and power,
and compare them with the test in \cite{CZZ2010a} for
testing $H_0\dvtx \Sigma=I_p$ and the test in \cite{CJ2011} for
testing a banded structure.

First we consider testing $H_0\dvtx  \Sigma=I_p$ against $H_a\dvtx  \Sigma\neq
I_p$. Draw $1000$ random samples with sample size $n=50$ or $200$
from the random variable $W_1+(\delta/n^{1/4})W_2$, where
$W_1\sim N(0,I_p)$, $W_2\sim N(0,(\sigma_{ij})_{1\le i, j\le p})$ with
$\sigma_{ij}=0.5^{|i-j|}\*I(|i-j|<\tau)$, and $W_1$ is independent of
$W_2$. When the sample size is small, it turns out that the size of
the proposed empirical likelihood ratio test is a bit large, and some
calibration is necessary. Here we propose the following bootstrap
calibration for the empirical likelihood ratio function $L_2(I_p)$ in
Theorem \ref{main1}.

For a given sample $\{\boldR_i^*(I_{p})\}_{i=1}^N$, we draw 300
bootstrap samples with size~$N$, say
$\{\tilde\boldR_i^{*(b)}(I_p)\}_{i=1}^N$ with $b=1,\ldots,300$.
Based on each bootstrap sample
$\{\tilde\boldR_i^{*(b)}(I_p)\}_{i=1}^N$, we compute the
bootstrapped empirical likelihood ratio function
\begin{eqnarray*}
L_2^{(b)}(I_p)&=&\sup\Biggl\{\prod
_{i=1}^N(Np_i)\dvtx
p_1\ge0,\ldots,p_N\ge0, \sum
_{i=1}^Np_i=1,
\\
&&\hspace*{67.1pt}\sum_{i=1}^Np_i
\tilde\boldR_i^{*(b)}(I_p)=\frac1N\sum
_{j=1}^N\boldR_j^*(I_p)
\Biggr\}.
\end{eqnarray*}
Then the bootstrap calibrated
empirical likelihood ratio test with level $\gamma$ will reject the null
hypothesis $H_0\dvtx  \Sigma=I_p$ whenever $-2\log L_2(I_p)$ is larger than
the $[300(1-\gamma)]$th largest value of $\{-2\log
L_2^{(b)}(I_p)\}_{b=1}^{300}$. More details on calibration for
empirical likelihood ratio test can be found in \cite{O2001}. We denote
the empirical likelihood ratio test based on $-2\log L_2(I_p)$, its
bootstrap calibrated version and the test in \cite{CZZ2010a} by
$\operatorname{EL}(\gamma)$, $\operatorname{BCEL}(\gamma)$ and $\operatorname{CZZ}(\gamma)$,
respectively, where $\gamma$ is the significance level.

\begin{table}
\tabcolsep=2pt
\caption{Sizes and powers are reported for the proposed
empirical likelihood method ($\operatorname{EL}(\gamma)$), its bootstrap calibrated
version ($\operatorname{BCEL}(\gamma)$) and the test in \cite{CZZ2010a} ($\operatorname{CZZ}(\gamma
)$) with significance level $\gamma=0.05$ for tesing $H_0\dvtx  \Sigma
=I_p$. We choose $\tau=10$}
\label{Tab1}
\begin{tabular*}{\tablewidth}{@{\extracolsep{\fill}}lcccccc@{}}
\hline
&\multicolumn{1}{c}{$\bolds{\operatorname{EL}(0.05)}$}
& \multicolumn{1}{c}{$\bolds{\operatorname{BCEL}(0.05)}$}
& \multicolumn{1}{c}{$\bolds{\operatorname{CZZ}(0.05)}$}
& \multicolumn{1}{c}{$\bolds{\operatorname{EL}(0.05)}$}
& \multicolumn{1}{c}{$\bolds{\operatorname{BCEL}(0.05)}$}
& \multicolumn{1}{c@{}}{$\bolds{\operatorname{CZZ}(0.05)}$}
\\
$\bolds{(n,p)}$ & \multicolumn{1}{c}{$\bolds{\delta=0}$}
& \multicolumn{1}{c}{$\bolds{\delta=0}$} & \multicolumn{1}{c}{$\bolds{\delta=0}$}
& \multicolumn{1}{c}{$\bolds{\delta=1}$} & \multicolumn{1}{c}{$\bolds{\delta=1}$}
& \multicolumn{1}{c@{}}{$\bolds{\delta=1}$}\\
\hline
$(50,25)$&0.127&0.054&0.053&0.296&0.118&0.219\\
$(50,50)$&0.148&0.065&0.067&0.324&0.136&0.216\\
$(50,100)$&0.138&0.068&0.038&0.317&0.125&0.212\\
$(50,200)$&0.168&0.081&0.041&0.310&0.113&0.221\\
$(50,400)$&0.151&0.071&0.045&0.342&0.145&0.242\\
$(50,800)$&0.154&0.064&0.041&0.337&0.137&0.219\\
[6pt]
$(200,25)$&0.065&0.048&0.052&0.348&0.305&0.179\\
$(200,50)$&0.058&0.052&0.041&0.336&0.298&0.162\\
$(200,100)$&0.068&0.054&0.059&0.353&0.319&0.179\\
$(200,200)$&0.056&0.051&0.058&0.358&0.322&0.155\\
$(200,400)$&0.069&0.064&0.051&0.374&0.343&0.180\\
$(200,800)$&0.058&0.047&0.050&0.366&0.338&0.182\\
\hline
\end{tabular*}
\end{table}

Table \ref{Tab1} reports the sizes ($\delta=0$) and powers ($\delta=1$) of
these three tests with level 0.05 by considering $\tau=10$ and $p=25,
50, 100, 200, 400, 800$. As we can see: (i) the empirical likelihood
ratio test
has a large size for the small sample size $n=50$, but the
bootstrap calibrated version has an accurate size, which is
comparable to the size of the test in \cite{CZZ2010a};
(ii) the test in \cite{CZZ2010a} is more powerful for
$n=50$, but less powerful when $n=200$; (iii) for a large sample
size, the empirical likelihood ratio test has no need to calibrate.

Next we consider testing $H_0\dvtx \sigma_{ij}=0$ for $|i-j|\ge\tau$ by
drawing $1000$ random samples from $\tilde W+(\delta/
n^{1/4})\bar W$, where $\tilde W\sim N(0,
(0.5^{|i-j|}I(|i-j|<\tau))_{1\le i,j\le p})$, $\bar
W=(\sum_{i=1}^{k}W_i/\sqrt k,\ldots,\sum_{i=p}^{p+k}W_i/\sqrt k)^T$,
$W_1,\ldots,W_{p+k}$ are i.i.d. with $N(0,1)$ and independent of
$\tilde W$. We consider the proposed empirical likelihood ratio test
based on
Theorem \ref{THCJb} ($\operatorname{EL}(\gamma)$) and a similar bootstrap calibrated
version as in testing $H_0\dvtx \Sigma=I_p$ ($\operatorname{BCEL}(\gamma)$), and compare
them with the test based on maximum in \cite{CJ2011}
($\operatorname{CJ}(\gamma)$).

Table \ref{Tab2} reports the sizes ($\delta=0$) and powers ($\delta=1)$ of
these three tests with level 0.05 by considering $\tau=5$,
$k=\tau+10$ and $p=25, 50, 100, 200, 400,\break 800$. From Table \ref{Tab2}, we
observe that: (i) the empirical likelihood ratio test has a large size for
the small sample size $n=50$, but the bootstrap calibrated version
has an accurate size, which is more accurate than the size of the
test in \cite{CJ2011}; (ii) the test in \cite{CJ2011}
has little power for all considered cases, and is much less powerful
than the proposed empirical likelihood ratio test; (iii) for a large sample
size, the empirical likelihood ratio test has no need to calibrate.

\begin{table}
\tabcolsep=2pt
\caption{Sizes and powers are reported for the proposed
empirical likelihood method ($\operatorname{EL}(\gamma)$), its bootstrap calibrated
version ($\operatorname{BCEL}(\gamma)$) and the test in \cite{CJ2011} ($\operatorname{CJ}(\gamma
)$) with significance level $\gamma=0.05$ for testing $H_0\dvtx  \sigma
_{ij}=0$ for all $|i-j|\ge\tau$. We choose $\tau=5$, $k=\tau+10$}
\label{Tab2}
\begin{tabular*}{\tablewidth}{@{\extracolsep{\fill}}lcccccc@{}}
\hline
&\multicolumn{1}{c}{$\bolds{\operatorname{EL}(0.05)}$}
& \multicolumn{1}{c}{$\bolds{\operatorname{BCEL}(0.05)}$}
& \multicolumn{1}{c}{$\bolds{\operatorname{CZZ}(0.05)}$}
& \multicolumn{1}{c}{$\bolds{\operatorname{EL}(0.05)}$}
& \multicolumn{1}{c}{$\bolds{\operatorname{BCEL}(0.05)}$}
& \multicolumn{1}{c@{}}{$\bolds{\operatorname{CZZ}(0.05)}$}
\\
$\bolds{(n,p)}$ & \multicolumn{1}{c}{$\bolds{\delta=0}$}
& \multicolumn{1}{c}{$\bolds{\delta=0}$} & \multicolumn{1}{c}{$\bolds{\delta=0}$}
& \multicolumn{1}{c}{$\bolds{\delta=1}$} & \multicolumn{1}{c}{$\bolds{\delta=1}$}
& \multicolumn{1}{c@{}}{$\bolds{\delta=1}$}\\
\hline
$(50,25)$&0.118&0.036&0.015&0.272&0.093&0.017\\
$(50,50)$&0.124&0.049&0.010&0.266&0.097&0.018\\
$(50,100)$&0.126&0.057&0.005&0.268&0.099&0.004\\
$(50,200)$&0.128&0.058&0.003&0.268&0.100&0.001\\
$(50,400)$&0.113&0.053&0.002&0.282&0.121&0.001\\
$(50,800)$&0.128&0.062&0.001&0.281&0.109&0.000\\
[6pt]
$(200,25)$&0.078&0.062&0.019&0.288&0.253&0.034\\
$(200,50)$&0.074&0.059&0.033&0.323&0.286&0.020\\
$(200,100)$&0.057&0.053&0.019&0.332&0.304&0.044\\
$(200,200)$&0.066&0.046&0.024&0.293&0.263&0.032\\
$(200,400)$&0.061&0.052&0.020&0.336&0.304&0.016\\
$(200,800)$&0.053&0.046&0.026&0.317&0.297&0.025\\
\hline
\end{tabular*}
\end{table}

It is expected that the test based on the maximum statistic in \cite
{CJ2011} should be more powerful than a test based on a Euclidean
distance when a large departure, instead of many small departures, from
the null hypothesis happens. To examine this, we test $H_0\dvtx  \sigma
_{ij}=0$ for $|i-j|\ge\tau$ by drawing 1000 random samples with size
$n=200$ from $\tilde W+\delta\bar W$, where $\tilde W\sim N(0,
(0.5^{|i-j|}I(|i-j|<\tau))_{1\le i,j\le p})$, $\bar W=(\bar
W_1,\ldots,\bar W_p)^T$ with $\bar W_1=\bar W_{\tau+1}\sim N(0,1)$ and
$\bar W_j=0$ for $j\neq1, \tau+1$. Again, $\tilde W$ and $\bar W_1$ are
independent. We take $\tau=5$, level 0.05 and $\delta=0.6, 0.7, 0.8$.
This is the sparse case in which we expect the $\mathrm{CJ}$ test to be favored.
The powers of $\operatorname{CJ}(0.05)$ are 0.074, 0.268 and 0.642 for $\delta=0.6,
0.7, 0.8$, respectively, while the powers of $\operatorname{EL}(0.05)$ are 0.066 for
all $\delta=0.6, 0.7, 0.8$. This confirms the advantage of using
maximum when a large departure occurs. However, as we argue in the
\hyperref[sec1]{Introduction}, the proposed empirical likelihood ratio test is quite
flexible in taking information into account. Since only one large
departure exists, the second equation in the proposed empirical
likelihood ratio test should be replaced by an estimating equation
related with this sparsity. Here, we use the first 40\% data to get the
sample variance $\hat\sigma_{ij}$ and find the positions of the largest
four values of $|\hat\sigma_{ij}|$ for $i-j\ge\tau$. Next we use the
remaining 60\% data to formulate the empirical likelihood ratio test
through replacing $\tilde v^{*\prime}$ in the second estimating
equation of $L_5(\Sigma)$ by the sum of values at the identified four
positions of the covariances $(\bold Y_i+\bold Y_{N+i})$. For this
modified empirical likelihood ratio test, we find that the empirical
size is 0.061, and powers are 0.106, 0.255 and 0.542 for $\delta=0.6,
0.7, 0.8$, respectively. As we can see, the empirical likelihood ratio
test with the new second equation improves the power significantly and
becomes comparable with the $\mathrm{CJ}$ test based on the maximum statistic.
In conclusion, the proposed empirical likelihood ratio test is powerful
and flexible.

\section{Proofs}\label{sec4}\label{proof}

Without loss of generality, we assume $\boldu_0=0$ throughout. For
simplicity, we use $\|\cdot\|$ to denote the $L_2$ norm of a vector or
matrix and write $\bolde_i(\Sigma_0)=\bolde_i$, $\boldv_i(\Sigma
_0)=\boldv_i$, $\bolde^*_i(\Sigma_0)=\bolde^*_i$, $\boldv
^*_i(\Sigma_0)=\boldv^*_i$,
$\bolde_i^{\prime}(\Sigma_0)=\bolde_i^{\prime}$, $\tilde\boldv
_i^{\prime}(\Sigma_0)=\tilde\boldv_i^{\prime}$, $\bolde
_i^{*\prime}(\Sigma_0)=\bolde_i^{*\prime}$, $\tilde\boldv
_i^{*\prime}(\Sigma_0)=\tilde\boldv_i^{*\prime}$,
$\pi_{11}=\E(\bolde_1^2(\Sigma_0))$ and $\pi_{22}=\E(\boldv
_1^2(\Sigma_0))$. We first collect some lemmas and leave the proofs in
the supplementary file.

\begin{lemma}\label{lem1}
Under condition (\ref{t21}) in Theorem \ref{main1}, we have
%
\begin{equation}
\label{411}{1\over\sqrt{N}}\sum_{i=1}^{N}
\biggl({\bolde_i\over
\sqrt{\pi_{11}}}, {\boldv_i\over
\sqrt{\pi_{22}}} \biggr)^T
\stackrel{d} {\longrightarrow} N(0, I_2).
\end{equation}
Further,
%
\begin{eqnarray}
\label{413} {{\sum_{i=1}^{N}\bolde
_i^2}\over{ N\pi_{11}}}-1 &\stackrel{p} {\longrightarrow}& 0,\qquad
{{\sum_{i=1}^{N}\boldv^2_i}\over{ N\pi_{22}}}-1 \stackrel{p}
{\longrightarrow} 0,\qquad
{{\sum_{i=1}^{N}\bolde_i\boldv_i}\over
N\sqrt{\pi_{11}\pi_{22}}} \stackrel{p} {\longrightarrow} 0,\hspace*{-28pt}
\\
\label{416}
\max_{1\leq i\leq
N}|\bolde_{i}/\sqrt{
\pi_{11}}|&=&o_p\bigl(N^{1/2}\bigr),\qquad  \max
_{1\leq
i\leq N}|\boldv_{i}/\sqrt{\pi_{22}}|=o_p
\bigl(N^{1/2}\bigr).
\end{eqnarray}
\end{lemma}

\begin{lemma}\label{lem1e*v*}
Under conditions (\ref{t21}) and (\ref{t21add}) in Theorem \ref
{main1}, we have
%
\begin{equation}
\label{412}{1\over
\sqrt{N}}\sum_{i=1}^{N}
\biggl({\bolde^*_i\over\sqrt{\pi_{11}}}, {\boldv^*_i\over
\sqrt{\pi_{22}}} \biggr)^T
\stackrel{d} {\longrightarrow} N(0, I_2).
\end{equation}
Further,
%
\begin{eqnarray}
\label{414} {{\sum_{i=1}^{N}\bolde
^{*2}_i}\over{N\pi_{11}}}-1 &\stackrel{p} {\longrightarrow}& 0,\qquad
{{\sum_{i=1}^{N}\boldv^{*2}_i}\over{N\pi_{22}}}-1 \stackrel{p}
{\longrightarrow} 0,\qquad
{{\sum_{i=1}^{N}\bolde^*_i\boldv^*_i}\over N\sqrt{\pi_{11}\pi_{22}}}
\stackrel{p} {\longrightarrow} 0,\hspace*{-28pt}
\\
\label{417}\max_{1\leq i\leq
N}\bigl|\bolde^*_{i}/\sqrt{
\pi_{11}}\bigr|&=&o_p\bigl(N^{1/2}\bigr),\qquad
\max_{1\leq
i\leq N}\bigl|\boldv^*_{i}/\sqrt{\pi_{22}}\bigr|=o_p
\bigl(N^{1/2}\bigr).
\end{eqnarray}
\end{lemma}

\begin{lemma}\label{lem2} Under conditions of Corollary \ref{main2},
for any $\delta>0$, we have
\[
\mathrm{E}|\bolde_1|^{2+\delta}\leq q^{\delta} \Biggl(\sum
_{i=1}^{p}\sum_{j=1}^p
\mathrm{E}|X_{1i}X_{1j}-\sigma_{ij}|^{2+\delta}
\Biggr)^2
\]
and
\[
\mathrm{E}|\boldv_1|^{2+\delta}\leq2^{4+\delta}q^{1+\delta}
\sum_{i=1}^{p}\sum
_{j=1}^p \mathrm{E}|X_{1i}X_{1j}-
\sigma_{ij}|^{2+\delta}.
\]
\end{lemma}

\begin{lemma}\label{lem3} Under conditions of Corollary \ref{main3},
we have
\[
\mathrm{E}\bolde_1^4/\bigl(\mathrm{E}
\bolde_1^2\bigr)^2=O(1) \quad\mbox{and}\quad
\mathrm{E}\boldv_1^4/\bigl(\mathrm{E}\boldv_1^2
\bigr)^2=O(1).
\]
\end{lemma}

\begin{lemma}\label{lem4} Under conditions of Theorem \ref{THCJb}, we
have
%
\begin{eqnarray}
\label{427} {\mathrm{E}\bolde_1'^4/\bigl(
\mathrm{E}\bolde_1'^2\bigr)^2}&=&O(1),\nonumber\\[-8pt]\\[-8pt]
{\mathrm{E}\tilde\boldv_1'^4/
\bigl(\mathrm{E}\tilde\boldv_1'^2
\bigr)^2}&=&O(1),\nonumber
\\
\label{l4add}\E\Biggl\{ \sum_{i=1}^{N}
\bigl(\bolde_i^{*\prime}-\bolde_i^{\prime}
\bigr)^2+ \Biggl[\sum_{i=1}^{N}
\bigl(\bolde_i^{*\prime}-\bolde_i^{\prime}
\bigr)\Biggr]^2\Biggr\} &=&o\bigl(N\mathrm{E}\bigl[
\bolde_1^{\prime2}\bigr]\bigr),
\\
\label{l4add2}\E\Biggl\{\sum_{i=1}^{N}\bigl(\tilde
\boldv_i^{*\prime}-\tilde\boldv_i^{\prime}
\bigr)^2+ \Biggl[\sum_{i=1}^{N}
\bigl(\tilde\boldv_i^{*\prime}-\tilde\boldv_i^{\prime}
\bigr)\Biggr]^2\Biggr\} &=& o\bigl(N\mathrm{E}\bigl[\tilde
\boldv_1^{\prime2}\bigr]\bigr).
\end{eqnarray}
\end{lemma}

\begin{pf*}{Proof of Theorem \ref{main1}}
Put $\hat e_i=e_i/\sqrt{\pi_{11}}$, $\hat v_i=v_i/\sqrt{\pi_{22}}$
and $\hat{\boldR}_i=(\hat e_i, \hat v_i)^T$ for $i=1,\ldots,N$. Then
it is easy to see that
$-2\log L_1(\Sigma_0)=2\times\sum_{i=1}^N\log\{1+\rho^T\hat{\boldR}_i\}$,
where $\rho=(\rho_1,\rho_2)^T$ satisfies
%
\begin{equation}
\label{pfmain1-1} \frac1N\sum_{i=1}^N
\frac{\hat{\boldR}_i}{1+\rho^T\hat{\boldR}_i}=0.
\end{equation}
Using Lemma \ref{lem1} and similar arguments in the proof of (2.14)
in \cite{O1990}, we can show that
%
\begin{equation}
\label{pfmain1-1a} \|\rho\|=O_p\bigl(N^{-1/2}\bigr).
\end{equation}
Then it follows from (\ref{416}) and (\ref{pfmain1-1a}) that
%
\begin{equation}
\label{pfmain1-1b} \max_{1\le i\le
N}\biggl|\frac{\rho^T\hat{\boldR}_i}{1+\rho^T\hat{\boldR}_i}\biggr|=o_p(1).
\end{equation}
By (\ref{pfmain1-1}), we have
\begin{eqnarray*}
0&=&\frac1N\sum_{i=1}^N\frac{\rho^T\hat{\boldR}_i}{1+\rho^T\hat
{\boldR}_i}
\\
&=&\frac1N\sum_{i=1}^N\rho^T
\hat{\boldR}_i\biggl\{1-\rho^T\hat{\boldR
}_i+\frac{(\rho^T\hat{\boldR}_i)^2}{1+\rho^T\hat{\boldR}_i}\biggr\}
\\
&=&\frac1N\sum_{i=1}^N\rho^T
\hat{\boldR}_i-\frac1N\sum_{i=1}^N
\bigl(\rho^T\hat{\boldR}_i\bigr)^2+\frac1N
\sum_{i=1}^N\frac{(\rho
^T\hat{\boldR}_i)^3}{1+\rho^T\hat{\boldR}_i}
\\
&=&\frac1N\sum_{i=1}^N\rho^T
\hat{\boldR}_i-\frac
{1+o_p(1)}{N}\sum_{i=1}^N
\bigl(\rho^T\hat{\boldR}_i\bigr)^2,
\end{eqnarray*}
which implies
%
\begin{equation}
\label{pfmain1-2} \frac1N\sum_{i=1}^N
\rho^T\hat{\boldR}_i=\frac{1+o_p(1)}N\sum
_{i=1}^N\bigl(\rho^T\hat{
\boldR}_i\bigr)^2.
\end{equation}
Using (\ref{pfmain1-1})--(\ref{pfmain1-1b}) and Lemma \ref{lem1}, we
have
\begin{eqnarray*}
0&=&\frac1N\sum_{i=1}^N\frac{\hat{\boldR}_i}{1+\rho^T\hat{\boldR
}_i}
\\[-2.5pt]
&=&\frac1N\sum_{i=1}^N\hat{
\boldR}_i\biggl\{1-\rho^T\hat{\boldR}_i+
\frac{(\rho^T\hat{\boldR}_i)^2}{1+\rho^T\hat{\boldR}_i}\biggr\}
\\[-2.5pt]
&=&\frac1N\sum_{i=1}^N\hat{
\boldR}_i-\frac1N\sum_{i=1}^N
\hat{\boldR}_i\hat{\boldR}_i^T\rho+\frac1N
\sum_{i=1}^N\frac{\hat
{\boldR}_i(\rho^T\hat{\boldR}_i)^2}{1+\rho^T\hat{\boldR}_i}
\\[-2.5pt]
&=&\frac1N\sum_{i=1}^N\hat{
\boldR}_i-\frac1N\sum_{i=1}^N
\hat{\boldR}_i\hat{\boldR}_i^T
\rho+O_p\Biggl(\max_{1\le i\le N}\biggl\|\frac{\hat
{\boldR}_i}{1+\rho^T\hat{\boldR}_i}\biggr\|
\frac1N\sum_{i=1}^N\bigl(\rho
^T\hat{\boldR}_i\bigr)^2\Biggr)
\\[-2.5pt]
&=&\frac1N\sum_{i=1}^N\hat{
\boldR}_i-\frac1N\sum_{i=1}^N
\hat{\boldR}_i\hat{\boldR}_i^T
\rho+o_p\Biggl(N^{1/2}\rho^T\frac1N\sum
_{i=1}^N\hat{\boldR}_i\hat{
\boldR}_i^T\rho\Biggr)
\\[-2.5pt]
&=&\frac1N\sum_{i=1}^N\hat{
\boldR}_i-\frac1N\sum_{i=1}^n
\hat{\boldR}_i\hat{\boldR}_i^T
\rho+o_p\bigl(N^{1/2}\bigr),
\end{eqnarray*}
which implies that
%
\begin{equation}
\label{pfmain1-3} \rho=\Biggl\{\frac1N\sum_{i=1}^N
\hat{\boldR}_i\hat{\boldR}_i^T\Biggr
\}^{-1}\frac1N\sum_{i=1}^N\hat{
\boldR}_i+o_p\bigl(N^{-1/2}\bigr).
\end{equation}
Hence, using Taylor expansion, (\ref{pfmain1-2}), (\ref{pfmain1-3})
and Lemma \ref{lem1}, we have
%
\begin{eqnarray}
\label{chi-sq} &&-2\log L_1(\Sigma_0)
\nonumber\hspace*{-25pt}
\\[-2.5pt]
&&\qquad=2\sum_{i=1}^N\rho^T\hat{
\boldR}_i-\bigl(1+o_p(1)\bigr)\sum
_{i=1}^N\bigl(\rho^T\hat{
\boldR}_i\bigr)^2
\nonumber\hspace*{-25pt}
\\[-2.5pt]
&&\qquad=\bigl(1+o_p(1)\bigr)\rho^T\sum
_{i=1}^N\hat{\boldR}_i\hat{
\boldR}_i^T\rho\hspace*{-25pt}
\\[-2.5pt]
&&\qquad=\bigl(1+o_p(1)\bigr) \Biggl(\frac1{\sqrt N}\sum
_{i=1}^N\!\hat{\boldR}_i
\Biggr)^{\!T}\!\Biggl(\frac1N\sum_{i=1}^N\!
\hat{\boldR}_i\hat{\boldR}_i^T
\Biggr)^{\!-1}\!\Biggl(\frac1{\sqrt{N}}\sum_{i=1}^N\!
\hat{\boldR}_i\Biggr)+o_p(1)
\nonumber\hspace*{-25pt}
\\[-2.5pt]
&&\qquad\stackrel{d} {\to}\chi^2_2 \qquad\mbox{as } n\to\infty.
\nonumber\hspace*{-25pt}
\end{eqnarray}
Similarly we can show that $-2\log
L_2(\Sigma_0)\stackrel{d}{\to}\chi^2_2$ by using Lemma
\ref{lem1e*v*}.\vadjust{\goodbreak}
\end{pf*}

\begin{pf*}{Proof of Corollary \ref{main2}}
First we prove the case of known $\mu$. Lem\-ma~\ref{lem2} implies
that under condition A2,
\[
\mathrm{E}|\bolde_1|^{2+\delta}=O\bigl(q^{2+\delta}\bigr)
\quad\mbox{and}\quad \mathrm{E}|\boldv_1|^{2+\delta}=O\bigl(q^{2+\delta}
\bigr).
\]
Further,
under condition A1, we have for a constant $C>0$,
$\pi_{11}=\operatorname{tr}(\Theta^2)\ge qC$ and $\pi_{22}={\mathbf
1}_{q}^{T}\Theta{\mathbf1}_q\ge qC$.
Thus,
\[
\mathrm{E}|\bolde_1|^{2+\delta}/\pi_{11}^{(2+\delta
)/2}=O
\bigl(q^{(2+\delta)/2}\bigr)=O\bigl(p^{2+\delta}\bigr)
\]
and
\[
\mathrm{E}|\boldv_1|^{2+\delta}/\pi_{22}^{(2+\delta
)/2}=O
\bigl(q^{(2+\delta)/2}\bigr)=O\bigl(p^{2+\delta}\bigr).
\]
Therefore, (\ref{t21}) in Theorem \ref{main1} follows from
condition A3, that is, Corollary~\ref{main2} holds for the case of
known $\mu$.

Next we prove the case of unknown $\mu$. Since (\ref{t21}) is
satisfied, by Theorem \ref{main1}, it is enough to show that condition
(\ref{t21add}) holds. Under condition $\max_{1\leq i\leq p}\sigma
_{ii}<C_0$, we have
%
\begin{equation}
\label{pt22a}\bigl(\operatorname{tr}\bigl(\Sigma^2\bigr)
\bigr)^2= \biggl(\sum_{1\le i, j\le p}
\sigma_{ij}^2 \biggr)^2\leq q^2
\Bigl(\max_{1\leq i\leq
p}\sigma_{ii}^2\Bigr)\leq
C_0^2q^2
\end{equation}
and
%
\begin{equation}
\label{pt22b} \bigl(\idv_p^T\Sigma\idv_p
\bigr)^2\leq q^2 \Bigl(\max_{1\leq i\leq p}
\sigma_{ii}^2\Bigr)\leq C_0^2q^2.
\end{equation}

On the other hand, under condition A1, there exists a constant $C>0$
such that
%
\begin{equation}
\label{pt22c}\pi_{11}=\operatorname{tr}\bigl(\Theta^2\bigr)\ge
qC \quad\mbox{and}\quad \pi_{22}={\mathbf1}_{q}^{T}
\Theta{\mathbf1}_q\ge qC.
\end{equation}
Note that condition A3 implies that $p=o(n^{1/4})$ and $q=o(n^{1/2})$.
Thus, by (\ref{pt22a}), (\ref{pt22b}) and (\ref{pt22c}), we have
\[
N\mathrm{E}e_1^2=N\pi_{11}\geq CNq\geq\bigl(
\operatorname{tr}\bigl(\Sigma^2\bigr)\bigr)^2
\]
and
\[
\sqrt{N}\mathrm{E}v_1^2=\sqrt{N}\pi_{22}\geq
\sqrt{N}qC> \bigl(\idv_p^T\Sigma\idv_p
\bigr)^2.
\]
Hence, (\ref{t21add}) holds and the proof of Corollary \ref{main2}
is complete.
\end{pf*}

\begin{pf*}{Proof of Corollary \ref{main3}}
It follows from
Lemma \ref{lem3} that (\ref{t21}) in Theorem \ref{main1} holds with
$\delta=2$. Hence Corollary \ref{main3} follows from Theorem~\ref{main1} when $\mu$ is known.

When $\mu$ is unknown, it follows from Lemma \ref{lem3} that (\ref
{t21}) holds. Further, through the proof of Lemma \ref{lem3}, we have
\[
\E\bigl[e_1^2\bigr]\ge C^2\bigl(\tr\bigl(
\Sigma^2\bigr)\bigr)^2 \quad\mbox{and}\quad \E
\bigl[v_1^2\bigr]\ge C\bigl(\mathbf{1}_p^T
\Sigma\mathbf{1}_p\bigr)^2,
\]
that is, condition (\ref{t21add}) holds. Thus, by Theorem \ref
{main1}, Corollary \ref{main3} holds for unknown $\mu$.
\end{pf*}

\begin{pf*}{Proof of Theorem \ref{banded}}
Since the required moment conditions are satisfied, it follows from the
same arguments as in the proof
of Theorem \ref{main1}.
\end{pf*}

\begin{pf*}{Proof of Theorem \ref{THCJb}}
Using Lemma \ref{lem4}, the proof of Theorem \ref{THCJb} follows from
the same arguments as in the proof
of Theorem \ref{main1}.
\end{pf*}

\begin{pf*}{Proof of Theorem \ref{power1}}
We only show the case of known $\mu$ since the case of unknown $\mu$
can be proved similarly.

First we consider the case of $\nu=o(N)$. Note that under the
alternative hypothesis $H_a$, $\mathrm{E}\boldY_1=\Sigma$ and write
for $1\leq i\leq N$,
\[
\bolde_i(\Sigma_0)=\bolde_i(\Sigma)+
\operatorname{tr}\bigl((\Sigma-\Sigma_0)^2\bigr)+ \operatorname{tr}
\bigl((\Sigma-\Sigma_0) (\boldY_i+\boldY_{N+i}-2
\Sigma)\bigr)
\]
and
$\boldv_i(\Sigma_0)=\boldv_i(\Sigma)+2{\mathbf1}_p^T(\Sigma-\Sigma
_0)\idv_p$,
where $q=p^2$. As a result, we have
%
\begin{eqnarray}
\label{powereq1} &&{1\over\sqrt{N}}\sum_{i=1}^{N}
\biggl({\bolde_i(\Sigma_0)\over
\sqrt{\pi_{11}}}, {\boldv_i(\Sigma_0)\over\sqrt{\pi_{22}}} \biggr)^T
\nonumber\\
&&\qquad={1\over\sqrt{N}}\sum_{i=1}^{N}
\biggl({\bolde_i(\Sigma)\over\sqrt{\pi_{11}}}, {\boldv_i(\Sigma)\over
\sqrt{\pi_{22}}} \biggr)^T+
\sqrt{N} (\zeta_{n1}, \zeta_{n2} )^T\\
&&\qquad\quad{}+
{1\over\sqrt{N}}\sum_{i=1}^{N} \bigl(
\eta_i(\Sigma), 0 \bigr)^T,\nonumber
\end{eqnarray}
where
$\eta_i(\Sigma)=\operatorname{tr}((\Sigma-\Sigma_0)(\boldY_i+\boldY
_{N+i}-2\Sigma))/\sqrt{\pi_{11}}$.
Since $
\mathrm{E}[\eta_i(\Sigma)]=0$ and
%
\begin{eqnarray}
\label{poweradd0} \mathrm{E}\bigl[\eta_i(\Sigma)
\bigr]^2&=&4\E\bigl(\operatorname{tr}\bigl((\Sigma-\Sigma_0) (
\boldY_1 -\Sigma)\bigr)^2\bigr)/\pi_{11}
\nonumber
\\
&\le&4\E\bigl(\operatorname{tr}\bigl((\Sigma-\Sigma_0)^2\bigr)
\operatorname{tr}\bigl((\boldY_1 -\Sigma)^2\bigr)\bigr)/
\pi_{11}
\\
&=&O\bigl[\operatorname{tr}\bigl((\Sigma-\Sigma_0)^2 \bigr)/
\sqrt{\pi_{11}}\bigr]=o(1),
\nonumber
\end{eqnarray}
we have
%
\begin{equation}\label{zetaon}\quad
\frac1N\sum_{i=1}^{N}\eta_i^2(
\Sigma)=o_p(1) \quad\mbox{and}\quad \frac{\max_{1\le i\le N}|\eta_i(\Sigma
)|}{\sqrt
N}\le\sqrt{
\frac{\sum_{i=1}^N\eta_i^2(\Sigma)}{N}}\stackrel{p}\to0.\hspace*{-28pt}
\end{equation}
Hence it follows from Lemma \ref{lem1} that
%
\begin{equation}
\label{023}V_N\stackrel{d} {\to} N(0, I_2),
\end{equation}
where
\[
V_N=\pmatrix{V_{N1}
\cr
V_{N2} }=\frac1{\sqrt
N}\sum_{i=1}^N \left\{\pmatrix{
\dfrac{\bolde_i({\Sigma_0})}{\sqrt{\pi_{11}}}
\vspace*{2pt}\cr
\dfrac{\boldv_i({\Sigma_0})}{\sqrt{\pi_{22}}} }-\pmatrix{\zeta_{n1}
\cr
\zeta_{n2}} \right\}.
\]
Put
$W_i=(\frac{\bolde_i({\Sigma_0})}{\sqrt{\pi_{11}}}, \frac{\boldv
_i({\Sigma_0})}{\sqrt{\pi_{22}}})^T$.
Then it follows from the proof of Theorem \ref{main1} that
%
\begin{eqnarray}
\label{pfth6-1} &&-2\log L_1(\Sigma_0)
\nonumber
\\
&&\qquad=\bigl(1+o_p(1)\bigr) \Biggl(\frac1{\sqrt N}\sum
_{i=1}^NW_i\Biggr)^T\Biggl(
\frac1N\sum_{i=1}^NW_iW_i^T
\Biggr)^{-1}\frac1{\sqrt N}\sum_{i=1}^NW_i+o_p(1)
\nonumber\\
&&\qquad=\bigl(1+\zeta_{n1}^2+\zeta_{n2}^2
\bigr)^{-1}\nonumber\\
&&\qquad\quad{}\times\bigl[\bigl(1+\zeta_{n2}^2\bigr)
(V_{N1}+\sqrt N\zeta_{n1})^2-2
\zeta_{n1}\zeta_{n2}(V_{N1}+\sqrt N
\zeta_{n1})
\\
&&\hspace*{43.7pt}\qquad\quad{} \times(V_{n2}+\sqrt N\zeta_{n2}) +\bigl(1+
\zeta_{n1}^2\bigr) (V_{n2}+\sqrt N
\zeta_{n2})^2\bigr]\nonumber\\
&&\qquad\quad{}+o_{p}(1)
\nonumber
\\
&&\qquad=(V_{N1}+\sqrt N\zeta_{n1})^2
\bigl(1+o_p(1)\bigr)+(V_{n2}+\sqrt N\zeta
_{n2})^2\bigl(1+o_p(1)\bigr)+o_p(1).
\nonumber
\end{eqnarray}

If the limit of $\nu=N(\zeta_{n1}^2+\zeta_{n2}^2)$, say $\nu_0$, is
finite, then it follows from (\ref{023}) and (\ref{pfth6-1}) that
$-2\log L_1(\Sigma_0)$ converges in distribution to a noncentral
chi-square distribution with two degrees of freedom and noncentrality
parameter $\nu_0$. If $\nu$ goes to infinite, the limit of the
right-hand side of (\ref{022}) is 1. By (\ref{pfth6-1}), we have
%
\begin{eqnarray}
\label{add1a}\quad &&-2\log L_1(\Sigma_0)
\nonumber
\\
&&\qquad\ge\biggl(\frac{N\zeta_{n1}^2}2-V_{N1}^2\biggr)
\bigl(1+o_p(1)\bigr)+\biggl(\frac{N\zeta
_{n2}^2}2-V_{N2}^2
\biggr) \bigl(1+o_p(1)\bigr)+o_p(1)
\\
&&\qquad=\frac{\nu
}{2}\bigl(1+o_p(1)\bigr)-\bigl(V_{N1}^2+V_{N2}^2
\bigr) \bigl(1+o_p(1)\bigr)+o_p(1)\stackrel{p} {\to}
\infty,
\nonumber
\end{eqnarray}
which implies that the limit of the left-hand side of (\ref{022}) is
also 1. Thus (\ref{022}) holds when $\nu=o(N)$.

For the case of $\liminf\nu/N>0$, we first consider the case of $\lim
\inf\zeta_{n2}^2>0$. Since ${\sum_{i=1}^Np_i\boldR_i(\Sigma_0)}=0$
implies that ${\sum_{i=1}^N p_i\nu_i(\Sigma_0)}=0$, we have
%
\begin{eqnarray}
\label{add2} L_1(\Sigma_0)&\le&\sup\Biggl\{\prod
_{i=1}^N(Np_i)\dvtx p_1\ge0,\ldots,p_N\ge0, \sum_{i=1}^Np_i=1,
\sum_{i=1}^Np_i
\boldv_i(\Sigma_0)=0\Biggr\}\hspace*{-16pt}
\nonumber\\[-4pt]\\[-12pt]
&=&\sup\Biggl\{\prod_{i=1}^N(Np_i)\dvtx p_1
\ge0,\ldots,p_N\ge0,\sum_{i=1}^Np_i=1,
\sum_{i=1}^Np_i
\frac{\boldv_i(\Sigma_0)}{\sqrt{\pi
_{22}}}=0\Biggr\}.\hspace*{-16pt}
\nonumber
\end{eqnarray}
Define
\begin{eqnarray*}
L^*(\theta)&=&\sup\Biggl\{\prod_{i=1}^N(Np_i)\dvtx
p_1\ge0,\ldots,p_N\ge0,\sum
_{i=1}^Np_i=1,\\
&&\hspace*{91pt}\sum
_{i=1}^Np_i\biggl(\frac{\boldv_i(\Sigma
_0)}{\sqrt{\pi_{22}}}-
\zeta_{n2}\biggr)=\theta\Biggr\}.
\end{eqnarray*}
Put $\theta^*=\frac{1}N\sum_{i=1}^N(\frac{\boldv_i(\Sigma
_0)}{\sqrt{\pi_{22}}}-\zeta_{n2})$. Then
%
\begin{equation}
\label{add30} {\log L^*\bigl(\theta^*\bigr)}=0.
\end{equation}
Since $\E\{\boldv_i(\Sigma_0)/\sqrt{\pi_{22}}-\zeta_{n2}\}=\E\{
\boldv_i(\Sigma)/\sqrt{\pi_{22}}\}=0$
and $\E\{\boldv_i(\Sigma_0)/\sqrt{\pi_{22}}-\zeta_{n2}\}^2=1$
under $H_a\dvtx  \Sigma\neq\Sigma_0$,
we have by using Chebyshev's inequality that
%
\begin{equation}
\label{add3}P\bigl(\bigl|\theta^*\bigr|>N^{-2/5}\bigr)\to0.
\end{equation}
Using $\E\{\boldv_i(\Sigma_0)/\sqrt{\pi_{22}}-\zeta_{n2}\}^2=1$,
similar to the proof of (\ref{add1a}), we can show that
\[
-2\log L^*\bigl(\theta^*_1\bigr)\stackrel{p}
{\to}\infty\quad\mbox{and}\quad
{-}2\log L^*\bigl(\theta^*_2\bigr)\stackrel{p} {\to}\infty,
\]
where $\theta^*_1=N^{-1/4}$ and $\theta^*_2=-N^{-1/4}$, which satisfy
$N(\theta^*_1)^2=o(N)$ and\break $N(\theta^*_2)^2=o(N)$. It follows from
\cite{HL1990} that the set $\{\theta\dvtx  -2\log L^*(\theta)\le c\}
=:I_c$ is convex for any $c$.
Take $c=\min\{-2\log L^*(\theta_1^*), -2\log L^*(\theta_2^*)\}/2$.
By (\ref{add30}), we have that $\theta^*\in I_c$. Thus, if $-\zeta
_{n2}\in I_c$, then $-a\zeta_{n2}+(1-a)\theta^*\in I_c$ for any $a\in
[0, 1]$, which implies that one of $\theta_1^*$ and $\theta_2^*$ must
belong to $I_c$. As a result,
we have
\begin{eqnarray*}
&&P\bigl(\bigl|\theta^*\bigr|\le N^{-2/5}, -\zeta_{n2}\in
I_c\bigr)\\
&&\qquad\le P\bigl(\theta_1^*\in I_c
\mbox{ or } \theta_2^*\in I_c\bigr)
\\
&&\qquad=P\bigl(\min\bigl\{-2\log L^*\bigl(\theta_1^*\bigr), -2\log L^*
\bigl(\theta_2^*\bigr)\bigr\}=0\bigr)\to0,
\end{eqnarray*}
which, together with (\ref{add3}), implies
\begin{eqnarray*}
&&
P\bigl( -2\log L^*(-\zeta_{n2})>c\bigr)
\\
&&\qquad = P(-\zeta_{n2}\notin I_c)\\
&&\qquad\ge1- P\bigl(\bigl|\theta^*\bigr|\le
N^{-2/5}, -\zeta_{n2}\in I_c\bigr)\\
&&\qquad\quad{}-P\bigl(\bigl|
\theta^*\bigr|> N^{-2/5}\bigr)\to1,
\end{eqnarray*}
and therefore
%
\begin{equation}
\label{add33} -2\log L^*(-\zeta_{n2})\stackrel{p} {\to}\infty
\end{equation}
since $c\stackrel{p}{\to}\infty$.
Hence, combining with (\ref{add2}), we have
\[
P\bigl(-2\log L_1(\Sigma_0)>\xi_{1-\alpha}\bigr)\ge
P\bigl(-2\log L^*(-\zeta_{n2})>\xi_{1-\alpha}\bigr)\to1,
\]
when $\lim\inf\zeta_{n2}^2>0$.\vadjust{\goodbreak}

Next we consider the case of $\lim\inf\zeta_{n1}>0$. Define
\[
\pi_{33}=\E\bigl\{\operatorname{tr}\bigl((\Sigma-\Sigma_0) (
\boldY_i+\boldY_{N+i}-2\Sigma)\bigr)\bigr\}^2
\quad\mbox{and}\quad \zeta_{n3}=\frac{\operatorname{tr}((\Sigma-\Sigma_0)^2)}{\sqrt{\pi_{11}+\pi_{33}}}.
\]
As before, we have
%
\begin{eqnarray}
\label{add4} L_1(\Sigma_0)
&\le&\sup\Biggl\{\prod_{i=1}^N(Np_i)\dvtx p_1
\ge0,\ldots,p_N\ge0,\nonumber\\
&&\qquad\quad\hspace*{2pt} \sum_{i=1}^Np_i=1,
\sum_{i=1}^Np_i
\bolde_i(\Sigma_0)=0\Biggr\}
\nonumber\\[-8pt]\\[-8pt]
&&\qquad=\sup\Biggl\{\prod_{i=1}^N(Np_i)\dvtx p_1
\ge0,\ldots,p_N\ge0,\nonumber\\
&&\qquad\quad\hspace*{23.8pt}
\sum_{i=1}^Np_i=1,\sum_{i=1}^Np_i
\frac{\bolde_i(\Sigma_0)}{\sqrt{\pi
_{11}+\pi_{33}}}=0\Biggr\}.
\nonumber
\end{eqnarray}
Define
\begin{eqnarray*}
L^{**}(\theta)&=&\sup\Biggl\{\prod_{i=1}^N(Np_i)\dvtx
p_1\ge0,\ldots,p_N\ge0,\\
&&\hspace*{24pt}\sum
_{i=1}^Np_i=1,
\sum_{i=1}^Np_i\biggl(\frac
{\bolde_i(\Sigma_0)}{\sqrt{\pi_{11}+\pi_{33}}}-\zeta_{n3}\biggr)=\theta
\Biggr\}.
\end{eqnarray*}
Since $\bolde_1(\Sigma)$ and $\operatorname{tr}((\Sigma-\Sigma_0)(\boldY
_1+\boldY_{N+1}-2\Sigma))$ are
two uncorrelated variables with zero means, we have
\[
\operatorname{Var}\bigl(\bolde_1(\Sigma)+\operatorname{tr}\bigl((\Sigma-\Sigma
_0) (\boldY_1+\boldY_{N+1}-2\Sigma)\bigr)
\bigr)=\pi_{11}+\pi_{33}.
\]
As we have shown in the proof of Lemma \ref{lem3},
$\E|\bolde_1(\Sigma)|^4=o(N\pi_{11}^2)$. Following the same lines
for estimating $\E(\boldv_1^4)$ in the end of the proof of Lemma \ref
{lem3}, we
have
\[
\E\bigl\{\operatorname{tr}\bigl((\Sigma-\Sigma_0) (\boldY_1+
\boldY_{N+1}-2\Sigma)\bigr)\bigr\}^4=O\bigl(
\pi_{33}^2\bigr). %
\]
Then it follows that
\begin{eqnarray*}
&&\E\bigl\{\bolde_1(\Sigma)+\operatorname{tr}\bigl((\Sigma-
\Sigma_0) (\boldY_1+\boldY_{N+1}-2\Sigma)
\bigr)\bigr\}^4
\\
&&\qquad\le 8\bigl(\E\bigl|\bolde_1(\Sigma)\bigr|^4+\E\bigl\{\operatorname{tr}
\bigl((\Sigma-\Sigma_0) (\boldY_1+\boldY_{N+1}-2
\Sigma)\bigr)\bigr\} ^4\bigr)
\\
&&\qquad=o\bigl(N(\pi_{11}+\pi_{33})^2\bigr).
\end{eqnarray*}
Write
\[
\frac{\bolde_i(\Sigma_0)}{\sqrt{\pi_{11}+\pi_{33}}}-\zeta_{n3}=\frac
{\bolde_i(\Sigma)+\operatorname{tr}((\Sigma-\Sigma_0)(\boldY
_i+\boldY_{N+i}-2\Sigma))}{\sqrt{\pi_{11}+\pi_{33}}}.
\]
Then we have
\begin{eqnarray*}
\E\biggl(\frac{\bolde_i(\Sigma_0)}{\sqrt{\pi_{11}+\pi_{33}}}-\zeta
_{n3}\biggr)^4&=&
\frac{\E(\bolde_i(\Sigma)+\operatorname{tr}((\Sigma-\Sigma
_0)(\boldY_i+\boldY_{N+i}-2\Sigma)))^4}{(\pi_{11}+\pi_{33})^2}\\
&=&o(N).
\end{eqnarray*}
This ensures the validity of Wilks's theorem for $-2\log L^{**}(0)$;
that is,
$-2\*\log L^{**}(0)$ converges in distribution to a chi-square
distribution with one degree of freedom.
Similar to the proof of (\ref{add1a}), we can show that
\[
-2\log L^{**}\bigl(\theta^*_1\bigr)\stackrel{p} {\to}
\infty\quad\mbox{and}\quad {-}2\log L^{**}\bigl(\theta^*_2\bigr)
\stackrel{p} {\to}\infty,
\]
where $\theta^*_1=N^{-1/4}$ and $\theta^*_2=-N^{-1/4}$, which satisfy
$N(\theta^*_1)^2=o(N)$ and\break $N(\theta^*_2)^2=o(N)$.

Put $\theta^{**}=\frac{1}N\sum_{i=1}^N(\frac{\bolde_i(\Sigma
_0)}{\sqrt{\pi_{11}+\pi_{33}}}-\zeta_{n3})$. Then
%
\begin{equation}
\label{add50} {\log L^{**}\bigl(\theta^{**}\bigr)}=0.
\end{equation}
Since
\begin{eqnarray*}
&&\E\bigl\{\bolde_i(\Sigma_0)/\sqrt{\pi_{11}+
\pi_{33}}-\zeta_{n3}\bigr\}
\\
&&\qquad=\E\biggl\{\frac{\bolde_i(\Sigma)+\operatorname{tr}((\Sigma-\Sigma
_0)(\boldY_i+\boldY_{N+i}-2\Sigma))}{\sqrt{\pi_{11}+\pi_{33}}}\biggr\}=0
\end{eqnarray*}
and
\[
\E\biggl\{\frac{\bolde_i(\Sigma_0)}{\sqrt{\pi_{11}+\pi_{33}}}-\zeta
_{n3}\biggr\}^2=\E
\biggl\{\frac{\bolde_i(\Sigma)+\operatorname{tr}((\Sigma-\Sigma
_0)(\boldY_i+\boldY_{N+i}-2\Sigma))}{\sqrt{\pi_{11}+\pi_{33}}}\biggr\} ^2=1
\]
under $H_a\dvtx  \Sigma\neq\Sigma_0$,
we have from Chebyshev's inequality that
%
\begin{equation}
\label{add5}P\bigl(\bigl|\theta^{**}\bigr|>N^{-2/5}\bigr)\to0.
\end{equation}
By (\ref{poweradd0}), we have
$
\pi_{33}/\pi_{11}=O(\zeta_{n1})$, which implies that there exists a
constant $M>0$ such that
\[
\zeta_{n3}/N^{-1/4}=N^{1/4}\zeta_{n1}
\frac{\sqrt{\pi_{11}}}{\sqrt{\pi_{11}+\pi_{33}}}\ge N^{1/4}\zeta_{n1}\{1+M
\zeta_{n1}\} ^{-1/2}\to\infty
\]
since $\lim\inf\zeta_{n1}>0$.

Using (\ref{add50}), (\ref{add5}) and the same arguments in proving
(\ref{add33}), we have
$-2\log L^{**}(-\zeta_{n3})\stackrel{p}{\to}\infty$.
Hence, combining with (\ref{add4}), we have
\[
P\bigl(-2\log L_1(\Sigma_0)>\xi_{1-\alpha}\bigr)\ge
P\bigl(-2\log L^{**}(-\zeta_{n3})>\xi_{1-\alpha}\bigr)
\to1,
\]
when $\lim\inf\zeta_{n1}^2>0$.
Therefore (\ref{022}) holds when $\lim\inf\zeta_{n1}>0$.
This completes the
proof of Theorem \ref{power1}.
\end{pf*}

\begin{pf*}{Proof of Theorem \ref{power2}}
The proof is similar to that of Theorem \ref{power1}.
\end{pf*}

\section*{Acknowledgments}

We thank the Editor Professor Runze Li, an Associate Editor and two
reviewers for their
constructive comments.

\begin{supplement}
\stitle{Supplement to ``Tests for covariance matrix with fixed or
divergent dimension''}
\slink[doi]{10.1214/13-AOS1136SUPP} 
\sdatatype{.pdf}
\sfilename{aos1136\_supp.pdf}
\sdescription{This supplementary file contains detailed proofs of
Lemmas \ref{lem1}--\ref{lem4} used in Section \ref{proof}.}
\end{supplement}


\printaddresses

\end{document}